\newfont{\Bbb}{msbm10 scaled\magstephalf}
\documentclass[12pt]{amsart}
\usepackage{amsmath, amsthm, amscd, amsfonts}
\newcommand{\RR}{{\mathbb R}}

\newcommand{\CC}{{\mathbb C}}

\newcommand{\NN}{{\mathbb N}}

\newcommand{\BB}{{\mathbb B}}
\newcommand{{\Li}}{{\mathcal L}_p}
\newcommand{{\Bl}}{{\mathcal B}_p}

\newtheorem{theorem}{Theorem}[section]

\newtheorem{corollary}[theorem]{Corollary}

\def\bege{\begin{equation}} \def\ende{\end{equation}}

\def\begr{\begin{eqnarray}} \def\endr{\end{eqnarray}}

\def\bnum{\begin{enumerate}}  \def\enum{\end{enumerate}}

\begin{document}

\title[Compact weighted composition operators]
{Compact weighted composition operators and fixed points in convex
domains}

\author[D. Clahane]{Dana D. Clahane}
\address{
\newline Department of Mathematics\newline
University of California
\newline Riverside, CA 92521.}

\email{dclahane@math.ucr.edu}
\keywords{weighted composition
operator, bounded operator, compact operator, convex domain, Hilbert
space, weighted Hardy space, Bergman space, holomorphic, several
complex variables}

\subjclass[2000]{Primary: 47B38; Secondary: 26A16, 32A16, 32A26,
32A30, 32A37, 32A38, 32H02, 47B33.}

\date{August 25, 2005}

\begin{abstract}
We extend a classical result of Caughran/H. Schwartz and another
recent result of Gunatillake by showing that if $D$ is a bounded,
convex domain in $\CC^n$, $\psi:D\rightarrow\CC$ is analytic and
bounded away from zero toward the boundary of $D$, and
$\phi:D\rightarrow D$ is a holomorphic map such that the weighted
composition operator $W_{\psi,\phi}$ is compact on a holomorphic,
separable, functional Hilbert space (containing the polynomial
functions densely) on $D$ with reproducing kernel $K$ satisfying
$K(z,z)\rightarrow\infty$ as $z\rightarrow\partial D$, then $\phi$
has a unique fixed point in $D$.  We apply this result by making a
reasonable conjecture about the spectrum of $W_{\psi,\phi}$ based on
previous results.
\end{abstract}

\maketitle

\section{Introduction}

Let $\phi$ be a holomorphic self-map of a bounded domain $D$ in
$\CC^n$, and suppose that $\psi$ is a holomorphic function on $D$.
We define the linear operator $W_{\psi,\phi}$ on the linear space of
complex-valued, holomorphic functions ${\mathcal H}(D) $ by
\[
W_{\psi,\phi}(f)=\psi(f\circ\phi).
\]
$W_{\psi,\phi}$ is called the {\em weighted composition operator}
induced by the {\em weight symbol} $\psi$ and {\em composition
symbol} $\phi$.  Note that $W_ {\psi,\phi}$ is the (unweighted)
composition operator $C_\phi$ given by $C_\phi(f)=f\circ\phi$ when
$\psi=1$.

Weighted composition operators are not only fundamental ojects of
study in analysis, but they have recently been the subject of an
increasing amount of attention due to the fact that a rather vast
amount of results have been obtained about unweighted composition
operators in one and several complex variables.  It should also be
pointed out that by the Banach-Stone Theorem, the isometries between
the spaces of continuous functions on Hausdorff spaces turn out to
be weighted composition operators \cite{s}.  This phenomenon also
occurs when these spaces are replaced by other function spaces on
domains in $\CC$ or $\CC^n$ (cf. \cite{Banach}, \cite{cwb},
\cite{dr}, \cite{fj}, \cite{fo}) and continues to be studied.

It is natural to consider the dynamics of the sequence of iterates
of a composition symbol of a weighted composition operator and the
spectra of such operators.  The following classical result of
Caughran/Schwartz \cite{Caughran/Schwartz} began this line of
investigation.
\begin{theorem}\label{one} Let $\phi:\Delta\rightarrow\Delta$ be an analytic self-map of the
unit disk $\Delta$ in $\CC$.  If $C_\phi$ is compact or
power-compact on the Hardy space $H^2(\Delta)$, then the following
statements hold:

(a) $\phi$ must have a unique fixed point in $\Delta$ ({\em this
point turns out to be the so-called Denjoy-Wolff point $a$ of $\phi$
in $\Delta$ (see \cite[~Ch.~2]{Cowen/MacCluer})}.

(b) The spectrum of $C_\phi$ is the set consisting of $0$, $1$, and
all powers of $\phi'(a)$.
\end{theorem}
The analogue of this result for Hardy spaces of the unit ball
$\BB_n$ in $\CC^n$ was obtained by MacCluer in \cite{MacCluer100}:
\begin{theorem}\label{two} Let $\phi:\BB_n\rightarrow\BB_n$ be a holomorphic self-map of $\BB_n$, and suppose that
$p\geq 1$. If $C_\phi$ is compact or power compact on the Hardy
space $H^p(\BB_n)$, then

(a) $\phi$ must have a unique fixed point in $\BB_n$ ({\em again,
this point is the so-called Denjoy-Wolff point $a$ of $\phi$ in
$\BB_n$ (see \cite[~Ch.~2]{Cowen/MacCluer})}.

(b) The spectrum of $C_\phi$ is the set consisting of $0$, $1$, and
all products of eigenvalues of $\phi'(a)$.
\end{theorem}
This result also holds for weighted Bergman spaces of $\BB_n$
\cite{Cowen/MacCluer}.  The proofs of parts (a) of Theorems
\ref{one} and \ref{two} appeal to the Denjoy-Wolff theorems in
$\Delta$ and $\BB_n$.  It is natural to consider whether Theorem
\ref{one} holds when $\BB_n$ is replaced by more general bounded
symmetric domains or even the polydisk $\Delta^n$.
 Chu/Mellon in \cite{Chu/Mellon} recently showed that the
Denjoy-Wolff theorem fails in $\Delta^n$ for $n>1$; nevertheless, it
is shown in \cite{cspec} that MacCluer's results can be generalized
from $\BB_n$ to arbitrary bounded symmetric domains that are either
reducible or irreducible.

G. Gunatillake in the forthcoming paper \cite{g} has extended
Theorem \ref{one} to weighted composition operators on a certain
class of weighted Hardy spaces of $\Delta$ when $\psi$ is bounded
away from $0$ toward the unit circle in $\CC$:
\begin{theorem}\label{three}
Let $(b_j)_{j\in\NN}$ be a sequence of positive numbers such that
$\liminf_{j\rightarrow\infty} b_j^{1/j}\geq 1$, and let
$H^2_b(\Delta)$ be the weighted Hardy space of analytic functions
$f:\Delta\rightarrow\CC$ whose MacClaurin series
$f(z)=\sum_{j=0}^\infty a_jz^j$ satisfy $\sum_{j=0}^\infty
|a_j|^2b_j^2<\infty$.  Suppose that $\phi:\Delta\rightarrow\Delta$
is analytic, and let $\psi:\Delta\rightarrow\CC$ be an analytic map
that is bounded away from zero toward the unit circle.  Assume that
$W_{\psi,\phi}$ is compact on $H^2_b(\Delta)$.  Then the following
statements are true:

(a) $\phi$ has a unique fixed point $a\in\Delta$.

(b) The spectrum of $W_{\psi,\phi}$ is the set
\[\{0,\psi(a)\}\cup\{\psi(a)[\phi'(a)]^j:j\in\NN\}.
\]
\end{theorem}
Related results for unweighted Hardy spaces and an example
$W_{\psi,\phi}$ that is compact even though $C_\phi$ is not appear
in \cite{Shapiro/Smith}.

In what follows, we will introduce some basic notation. Then, we
will make a conjecture concerning a multivariable analogue of the
above result.  The purpose of this paper is to report that part (a)
of Theorem \ref{three} extends quite generally to a wide class of
functional Hilbert spaces on convex domains in one or more
variables.  The calculation of the spectrum of a compact
$W_{\psi,\phi}$ in the multivariable setting is the subject of
future work.

Some of the following notation and definitions are standard, but we
include them for the sake of clarity and completeness:

\section{Notation and definitions}

Fix $n\in\NN$.  We denote the usual Euclidean distance from
$z\in\CC^n$ to $A\subset\CC^n$ by $d(z,A)$, and we say that
$z\rightarrow A$ if and only if $d(z,A)\rightarrow 0$. Let $D$ be a
bounded domain in $\CC^n$. Suppose that $\nu\in\RR$, and assume that
$u:D\rightarrow \RR$. We write
\[
\liminf_{z\rightarrow\partial D}u(z)=\nu
\]
if and only if for all $\varepsilon>0$ there is a $\delta>0$ such
that whenever $0<d(z,\partial D)<\delta$, we have
\[
\left|\inf_{\{z\in D:d(z,\partial
D)<\delta\}}u(z)-\nu\right|<\varepsilon.
\]
Let $\psi:D\rightarrow\CC$.  We say that {\em $\psi$ is bounded away
from zero toward the boundary of $D$} iff there is a $\nu>0$ such
that
\[
\liminf_{z\rightarrow\partial D}|\psi(z)|=\nu.
\]
Given a functional Hilbert space ${\mathcal Y}$ of holomorphic
functions defined on a domain $D$, in $\CC^n$, the Riesz
representation theorem guarantees that for each $z\in D$, there is a
unique $K_z\in{\mathcal Y}$ such that
\[
f(z)=\langle f,K_z\rangle\,\,\,\,\text{ for all }f\in{\mathcal Y}.
\]
This uniqueness ensures allows to us to define the {\em reproducing
kernel} $K:D\times D\rightarrow\CC$ for ${\mathcal Y}$ by
$K(z,w)=K_z(w)$.

\section{A conjecture and the result}

Based on the results to date, we propose the natural problem of
resolving whether or not the following statement holds:

\medskip

{\bf Conjecture:} {\em Suppose that $D\subset\CC^n$ is a bounded,
convex domain such that a given functional Hilbert space of
holomorphic functions ${\mathcal Y}$ in which the polynomials are
contained densely has reproducing kernel $K$ satisfying
$K(z,z)\rightarrow \infty$ as $z\rightarrow\partial D$. Let
$\psi:D\rightarrow\CC$ be holomorphic and suppose that $\psi$ is
bounded away from $0$ toward $\partial D$.  Assume that
$\phi:D\rightarrow D$ is a holomorphic map and that $W_{\psi,\phi}$
is compact on ${\mathcal Y}$.  Then

(1) $\phi$ has a unique fixed point in $D$, and

(2) the spectrum of $W_{\psi,\phi}$ is the set
$\{\psi(a)\sigma:\sigma\in E\}$, where $E$ is the set consisting of
$0$, $1$, and all possible products of eigenvalues of $\phi'(a)$.}

\medskip

While (2) has not yet been resolved in the multivariable case, we
are able to report that item (1) indeed holds.  What is unique about
this result is that most results about properties of composition
operators have heavily depended on the function space under
consideration; contrastingly, by proving (1), we obtain a result
that is relatively independent of the given Hilbert space.

The following theorem generalizes results in
\cite{Caughran/Schwartz}, \cite{g}, and, in one direction, the fixed
point portion of \cite[~Thm.~4.2]{cspec}. The reader is referred to
\cite{cspec} for the definition of compact divergence, which is used
in the proof that follows.
\begin{theorem}\label{fixed point}
Let $D\subset\CC^n$ be a convex domain, and suppose that ${\mathcal
Y}$ is a functional Hilbert space of functions on $D$ with
reproducing kernel $K:D\times D\rightarrow\CC$.  Assume that
$K(z,z)\rightarrow\infty$ as $z\rightarrow\partial D$, and assume
that the polynomial functions on $D$ are dense in ${\mathcal Y}$.
Suppose that $\psi:D\rightarrow\CC$ is holomorphic and bounded away
from zero toward the boundary of $D$, and let $\phi:D\rightarrow D$
be holomorphic.  Assume that $W_{\psi,\phi}$ is compact on
${\mathcal Y}$.  Then $\phi$ has a unique fixed point in $D$.
\end{theorem}
\begin{proof}
Let $k_z=K_z/||K_z||_{\mathcal Y}$.  Since $K(z,z)\rightarrow\infty$
as $z\rightarrow\partial D$ and the polynomials functions on $D$ are
dense in ${\mathcal Y}$, one can show using an argument identical to
that of the proof of \cite[~Lemma~3.1]{cspec} that $k_z\rightarrow
0$ weakly as $z\rightarrow\partial D$.  From the linearity of
$W_{\psi,\phi}$ and the identity
\[
W_{\psi,\phi}^*K_z={\overline{\psi(z)}}K_{\phi(z)},
\]
it immediately follows that
\[
||W_{\psi,\phi}^*k_z||_{\mathcal
Y}^2=|\psi(z)|^2K(z,z)^{-1}K[\phi(z),\phi(z)].
\]
Since $k_z\rightarrow 0$ weakly as $z\rightarrow\partial D$, we then
have that
\begin{equation}\label{genzhu}
\lim_{z\rightarrow \partial
D}|\psi(z)|^2K(z,z)^{-1}K[\phi(z),\phi(z)]=0.
\end{equation}
First suppose that $\phi$ has no fixed point in $D$.  We will obtain
a contradiction.  Let $z\in D$.  Since $D$ is convex, the sequence
of iterates $\phi^{(j)}$ of $\phi$ is compactly divergent
\cite[~p.~274]{Abate}.  Thus, for every compact $K\subset D$ there
is an $N\in\NN$ such that $\phi^{(j)}(z)\in D\setminus K$ for all $
j\geq N$.

Since for any $\varepsilon>0$, the set $K_\varepsilon$ of all $w\in
D$ such that $d(w,\partial D)\geq\varepsilon$ is compact, it follows
from the statement above that for all $\varepsilon>0$ there is an
$N\in\NN$ such that for all $j\geq N$, $\phi^{(j)}(z)\notin
K_\varepsilon$; alternatively, $d(\phi^{(j)}(z),\partial
D)<\varepsilon$ for $j\geq N$. Hence, we have that
$\phi^{(j)}(z)\rightarrow\partial D$ as $j\rightarrow\infty$ for all
$ z\in D$.  Since $K(z,z)\rightarrow\infty$ as $z\rightarrow\infty$
by assumption, it must be the case that
\[
\lim_{j\rightarrow\infty}K[\phi^{(j)}(z),\phi^{(j)}(z)]=\infty.
\]
Consequently, for any $z\in D$, and for infinitely many values of
$j$, we have that
\begin{equation}\label{mag}
K\{\phi[\phi^{(j)}(z)],\phi[\phi^{(j)}(z)]\}>K[\phi^{(j)}(z),\phi^{(j)}(z)]>0.
\end{equation}
It follows from the assumption that $\psi$ is bounded away from $0$
toward the boundary of $D$ that there must be a $\delta>0$ such that
whenever $d(\xi,\partial D)<\delta$, $|\psi(\xi)|> \nu/2$.  In
addition, for sufficiently large $j$, we have that
$d(\phi^{(j)}(z),\partial D)<\delta$, so that for these values of
$j$, $|\psi[\phi^{(j)}(z)]|>\nu/2$.  Therefore, for any $z\in D$,
there is an $N\in\NN$ such that for infinitely many values of $j\geq
N$, the following inequality holds for infinitely many $j$'s:
\[
|\psi(\phi^{(j)}(z))|^2K\{\phi[\phi^{(j)}(z)],\phi[\phi^{(j)}(z)]\}>\frac{\nu^2}{4}K[\phi^{(j)}(z),\phi^{(j)}(z)]>0.
\]
Thus we have in particular that for any $z\in D$, there are
infinitely many values of $j$ such that
\[
|\psi[\phi^{(j)}(z)]|^2K[\phi^{(j)}(z),\phi^{(j)}(z)]^{-1}K\{\phi[\phi^{(j)}(z)],\phi[\phi^{(j)}(z)]\}>\frac{\nu^2}{4}.
\]
Denote this sequence of values of $j$ by $(j_k)_{k\in\NN}$.
 Then we have that $\phi^{(j_k)}(z)\rightarrow\partial D$ as $k\rightarrow\infty$.  This
 fact, in combination with the fact that the above inequality for
 any $w\in D$ holds for the subsequence $(j_k)_{k\in\NN}$ of $\NN$,
 leads to a contradiction of Equation $(\ref{genzhu})$.  Hence, the assumption that
 $\phi$ has no fixed points is false.

 To show that $\phi$ has only
 one fixed point, assume to the contrary that $\phi$ has more than
 one fixed point.  By a result of Vigu\'{e}, the fixed point set of a holomorphic self-map of a bounded, convex domain in $\CC^n$ is a connected,
 analytic submanifold of that domain (see \cite[~Thm.~4.1]{cspec} or \cite{Vigue}).  Note that since the fixed point set of $\phi$ is not a singleton, then
 this set is, in particular, uncountable.  Denote this set of fixed points by ${\mathcal F}$.  We then have
that
\begin{equation}\label{eigen}
W_{\psi,\phi}^*(K_a)={\overline{\psi(a)}}K_{\phi(a)}={\overline{\psi(a)}}K_a
\,\,\,{\text{for all}}\,\,\,a\in{\mathcal F}.
\end{equation}
Therefore, for all $a\in {\mathcal F}$, we have that
${\overline{\psi(a)}} $ is an eigenvalue of the compact operator
$W_{\psi,\phi}^*$.  Since $\psi$ is continuous and ${\mathcal F}$ is
a connected analytic manifold in $\CC^n$, it must be the case that
$\psi({\mathcal F})$ must be either a singleton or uncountable.

First, assume that $\psi({\mathcal F})$ is a singleton
$\{\lambda\}$, so that Condition (\ref{eigen}) becomes
\[
W_{\psi,\phi}^*(K_a)={\overline{\lambda}}K_{\phi(a)}={\overline{\lambda}}K_a
\,\,\,{\text{for all}}\,\,\,a\in{\mathcal F}.
\]
Suppose that $\lambda=0$.  Then ker$W_{\psi,\phi}^*$ has uncountable
dimension, since $\{K_a:a\in D\}$ is linearly independent, thus
contradicting the separability of ${\mathcal Y}$. Therefore, assume
that $\lambda\not=0$.  Since $\{K_a:a\in D\}$ is a linearly
independent set, it follows that the ${\overline{
\lambda}}$-eigenspace of $W_{ \psi,\phi}^*$ has infinite dimension.
However, by \cite[~Thm.~7.1]{Conway}, this infiniteness contradicts
the compactness of $W_{\psi,\phi}^*$ on ${\mathcal Y}^*$.

Next, assume that $\psi({\mathcal F})$ is uncountable.  Then by
Condition (\ref{eigen}), $W_{\psi,\phi}^*$ has uncountably many
eigenvalues ${\overline{\psi(a)}}$ with $a\in {\mathcal F}$. Now
since ${\mathcal Y}$ contains the polynomials and is, therefore,
infinite-dimensional, ${\mathcal Y}^*$ is also infinite-dimensional.
Therefore, the compact operator $W_{\psi,\phi}^*$ has a countably
many eigenvalues \cite[~p.~214]{Conway}, and we have again obtained
a contradiction.

Hence, it must be the case that our assumption that $\phi$ has more
than one fixed point is false.
\end{proof}

\section{Some remarks and related, open uestions}

(1) Note that if $D=\Delta$ or $\BB_n$, the fixed point of $\phi$ in
Theorem \ref{fixed point} is precisely the so-called Denjoy-Wolff
$a$ point of $\phi$, to which the iterates of $\phi$ converge
uniformly on compacta.  It is natural to extend this uniform
convergence to the general setting of Theorem \ref{fixed point}.

\medskip

(2) As stated in \cite{cspec}, an interesting aspect of the above
result is that in the case when $D=\Delta^n$, the Denjoy-Wolff
theorem fails, and there is no unique ``Denjoy-Wolff point".
Nevertheless, the above fixed point theorem holds even for reducible
convex domains such as $\Delta^n$.

\medskip

(2) The convexity of $D$ in the proof of Theorem \ref{fixed point}
is used in two places: (a) to establish that if $\phi$ has no
interior fixed points, the iterates of $\phi$ diverge compactly, and
(b) to establish the assertion that the fixed point set of $\phi$
when $D$ is convex is a {\em connected} analytic submanifold of $D$.
It is therefore of interest to determine to what extent the
hypothesis of convexity can be weakened in such a way that tasks (a)
and (b) can still be simultaneously completed.

\medskip

(3) Let $G$ be a simply connected region that is properly contained
in $\CC$, and suppose that $\tau:\Delta\rightarrow\CC$ is the
Riemann mapping for $G$.  Let $H^2(G)$ be the Hardy space of
functions $f:G\rightarrow\CC$ that are analytic and satisfy
\[
\sup_{0<r<1}\int_{\tau(\{z\in\Delta:|z|=r\})}|f(z)|^2|dz|<\infty.
\]
In \cite{Shapiro/Smith} it is shown that if $C_\phi$ is compact on
$H^2(G)$ for some analytic $\phi:\Delta\rightarrow\Delta$, then
$\phi$ must have a unique fixed point in $G$.  Of course, such a
domain $G$ can have boundary portions that are concave, though all
domains in $\CC$ are trivially pseudoconvex \cite{kbook}.  On the
other hand, as is well known, the Riemann mapping theorem does not
extend to several complex variables, and the proof in
\cite{Shapiro/Smith} does seem to rely on the Denjoy-Wolff theory
that is inherent from the convexity of $\Delta$.

\medskip

(4) Concerning the goal of obtaining compact divergence results for
domains that are not necessarily convex, we list some of the best
known results: First, in \cite{Huang}, it is proven by X. Huang that
the iterates of a holomorphic self-map of a topologically
contractible, strictly pseudoconvex domain form a compactly
divergent sequence; however, in \cite{a1}, M. Abate showed the
existence of a holomorphic self-map of a topologically contractible
{\em pseudoconvex} domain such that the map's iterates do not
compactly diverge.

\medskip

(5) Note that in the proof of Theorem \ref{fixed point}, all that
was needed from Vigu\'{e}'s theorem is the assertion that if the
fixed point set of a holomorphic self-map of a convex domain is
non-empty, then it either contains one point or uncountably many
points. Vigu\'{e} in \cite{v2} has shown that the fixed point set of
a holomorphic self-map of any bounded domain $D$ (note that
``convex" is omitted!) in $\CC^n$ is also an analytic submanifold of
$D$, but it is an interesting and open question as to whether or not
the fixed point set in this case is necessarily connected for
general bounded domains besides the convex ones.

M. Abate has conjectured that the answer is affirmative for a
topologically contractible, strictly pseudoconvex domain.  A
resolution of this conjecture, together with Huang's previously
mentioned result \cite{Huang} concerning compact divergence of
iterates of holomorphic self-maps on topologically contractible
pseudoconex domains, would imply that Theorem \ref{fixed point}
extends to these domains.

(6) It is natural to wonder if a holomorphic self-map of a domain in
$\CC^n$ can have more than one but finitely many fixed points. This
question is easy to answer if we allow the domain to be
non-contractible. Let $r\in(0,1)$, and let $D$ be the $n$-fold
topological product of the annulus with inner radius $r$ and outer
radius $1/r$.  Define $\phi:D\rightarrow D$ be given by
$\phi(z_1,z_2,\ldots,z_n)=(z_1^{-1},z_2^{-1},\ldots,z_n^{-1})$.
$\phi$ has exactly $2^n$ fixed points, which are precisely the
points whose entries are either $1$ or $-1$.

Thus, one can ask the following question: If $D$ is a domain in
$\CC^n$ and $k\in\NN$, $k\not=1$ is given, is there a holomorphic
self-map of $D$ with precisely $k$ fixed points?  In this direction,
a result of J.-P. Vigu\'{e} \cite[~Thm.~1.7.6]{Krantz} states that
an analytic self-map of a bounded domain in $\CC$ with three
distinct fixed points is the identity mapping.  The multivariable
situation is quite different in this respect.

It would also be interesting to know how structured the fixed point
set is in this case; for example, one can ask what are the possible
Euclidean or fractal dimensions of such a fixed point set in
pseudoconvex domains.

(7) For the weighted Hardy spaces $H^2_b(\Delta)$ of the unit disk
in $\Delta\in\CC$, the Hardy spaces $H^2(D)$, and weighted Bergman
spaces $A^2_\alpha(D)$, where $D$ is either $\BB_n$, $\Delta^n$, or
more generally, any bounded symmetric domain in its Harish-Chandra
realization (see \cite{cspec}), the reproducing kernel $K$ satisfies
$K(z,z)\rightarrow\infty $ as $z\rightarrow\Delta$ (respectively,
$z\rightarrow D$), so the following fact, which extends the fixed
point results in \cite{Caughran/Schwartz} and \cite{g}, is an
immediate consequence of Theorem \ref{fixed point}:
\begin{corollary}
Suppose that ${\mathcal Y}$ is either the Hardy space $H^2(D)$ or
the weighted Bergman space $A^2_\alpha(D)$ of a bounded symmetric
domain $D$ with $\alpha<\alpha_D$, where $\alpha_D$ is a certain
critical value that depends on $D$ {\em (cf. \cite{cspec})}, and
assume that $\psi:D\rightarrow\CC$ is analytic and bounded away from
zero. Suppose that $\phi:D\rightarrow D$ is holomorphic, and let
$W_{\psi,\phi}$ be compact on ${\mathcal Y}$. Then $\phi$ has a
unique fixed-point in $D$.  This result also holds when $D=\Delta$
and ${\mathcal Y}=H^2_b(\Delta)$.
\end{corollary}
\begin{proof}
The assertions about $H^2(D)$ and $A^2_\alpha(D)$ immediately follow
from Theorem \ref{fixed point} and the fact that their reproducing
kernels approach infinity along the diagonal $\{(z,z):z\in D\}$ as
$z\rightarrow D$ (see \cite{cspec}). The assertion about ${\mathcal
Y}=H^2_b(\Delta)$ also immediately follows from Theorem \ref{fixed
point} and the fact that the assumed condition on the sequence
$(b_j)_{j\in\NN}$ implies that the reproducing kernel $K$ for
$H^2_b(\Delta)$ satisfies the same singularity property toward the
boundary along the diagonal (cf. \cite{cspec}).
\end{proof}

\section{Acknowledgment}

The author would like to thank M. Abate for correspondences that led
to some portions of Remarks (4)-(6) above.  Thanks are also extended
to W. Sheng and T. Oikhberg for helpful comments.


\begin{thebibliography}{99}
\bibitem[Ab]{Abate} M. Abate.  {\em Iteration Theory of Holomorphic Maps on Taut Manifolds},
Mediterranean Press, Rende, 1989.
\bibitem[Ab2]{a1} M. Abate.  {\em Holomorphic actions on
contractible domains without fixed points}, Math. Z. {\bf 211}
(1992), 547-555.
\bibitem[B]{Banach} S. Banach. {\em Theorie des
Operations Lineares}, Chelsea, Warsaw, 1932.
\bibitem[CaSc]{Caughran/Schwartz} J. Caughran/H. J. Schwartz.  {\em Spectra of compact composition operators}, Proc. Amer.
Math. Soc. {\bf 51}, (1975), 127-130.
\bibitem[ChMe]{Chu/Mellon} C.-H. Chu/P. Mellon.  {\em Iteration of compact holomorphic maps on a Hilbert ball}, Proc.
Amer. Math. Soc. {\bf 125} (1997), no. 6, 1771-1777.
\bibitem[CiW]{cwb} J. A. Cima, W. R. Wogen. {\em On isometries of the Bloch space}, Illinois J. Math. {\bf 24}
(2) (1980), 313-316.
\bibitem[Cl2]{cspec} D. D. Clahane. {\em Spectra of compact composition
operators over bounded symmetric domains}, Integral Equations
Operator Theory {\bf 51} (2005), no. 1, 41-56.
\bibitem[Con]{Conway} J. B. Conway. {\em A Course in Functional
Analysis}, second edition, Springer-Verlag, New York, 1990.
\bibitem[CowMac]{Cowen/MacCluer} C. C. Cowen/B. D. MacCluer.  {\em Composition
Operators on Spaces of Analytic Functions}, CRC Press, Boca Roton,
1995.
\bibitem[dRW]{dr} K. DeLeeuw, W. Rudin, J. Wermer. {\em The isometries of some function
spaces}, Proc. Amer. Math. Soc. {\bf 11} (1960), 485-694.
\bibitem[FlJ]{fj} R. J. Fleming/ J. E. Jamison. {\em Isometries on Banach spaces: function spaces},
Monographs and Surveys in Pure and Applied Mathematics, vol. 129,
Chapman and Hall/CRC, London, Boca Raton, FL, 2002.
\bibitem[Fo]{fo} F. Forelli. {\em The isometries of $H^p$}, Canad. J. Math. {\bf 16} (1964), 721-728.
\bibitem[G]{g} G. Gunatillake.  {\em The spectrum of a compact
weighted composition operator}, Proc. Amer. Math. Soc., to appear.
\bibitem[H]{Huang} X. Huang.  {\em A non-degeneracy property of
extremal mappings and iterates of holomorphic self-mappings}, Ann.
Scuola Norm. Pisa Cl. Sci. (4) {\bf 21} (1994), no. 3, 399-413.
\bibitem[K1]{Krantz} S. G. Krantz. {\em Geometric Analysis and Function
Spaces}, CBMS Regional Conference Series in Mathematics, {\bf 81},
Amer. Math. Soc., Providence, 1993.
\bibitem[K2]{kbook} S. G. Krantz.  {\em Function Theory of Several
Complex Variables}, Amer. Math. Soc., Providence, 2001.
\bibitem[Mac2]{MacCluer100} B. D. MacCluer.  {\em Spectra of compact composition
operators on $H^p(B_n)$}, Analysis {\bf 4} (1984), 87-103.
\bibitem[ShSm]{Shapiro/Smith} J. H. Shapiro/W. Smith. {\em Hardy spaces that support no compact composition
operators}, J. Funct. Anal. {\bf 205} (2003), no. 1, 62-89.
\bibitem[S]{s} M. Stone. {\em Application of the theory of Boolean rings in topology},
Trans. Amer. Math. Soc. {\bf 41} (1937) 375-481.
\bibitem[Vi1]{v2} J.-P. Vigu\'{e}. {\em Sur les points fixes d'applications holomorphes. (French) [On the fixed points
of holomorphic mappings]}, C. R. Acad. Sci. Paris Sér. I Math. {\bf
303}, (1986), no. 18, 927-930.
\bibitem[Vi2]{Vigue} J.-P. Vigu\'{e}.  {\em Points fixes d'applications holomorphes dans un domaine borné convexe de
$C\sp n$} [Fixed points of holomorphic mappings in a bounded convex
domain in ${\bf C}\sp n$], Trans. Amer. Math. Soc. {\bf 289} (1985),
no. 1, 345-353.
\end{thebibliography}
\end{document}